\documentclass[stropt,envcountsect,envcountsame, 12pt]{svjour3}
\usepackage[centertags]{amsmath}
\usepackage{graphicx}
\usepackage{amsfonts}
\usepackage{amssymb}
\usepackage{mathrsfs}
\usepackage{enumerate}
\usepackage[T1]{fontenc}

\newcommand{\ds}{\displaystyle}
\newcommand{\R}{\mathbb R}

\newcommand{\norm}[1]{\left\Vert#1\right\Vert}
\newcommand{\normsb}[1]{\norm{#1}}

\newcommand{\set}[2]{\left\{ #1 \,;\, #2 \right\}}
\newcommand{\eps}{\varepsilon}

\newcommand{\PS}[2]{\langle  #1 \mid #2  \rangle}

\newcommand{\vertiii}[1]{{\left\vert\kern-0.25ex\left\vert\kern-0.25ex\left\vert #1 
    \right\vert\kern-0.25ex\right\vert\kern-0.25ex\right\vert}}

\def\grad{\nabla}
\def\x{x}

\def\H{H}

\def\C{{\mathscr C}}
\def\CC{{\mathscr C}}

\def\solf{u}
\def\solfz{u_0}
\def\solV{v}

\def\v{v}
\def\solgen{v}

\def\xs{\hat{x}}
\def\xsp{\hat{x}_\star}

\def\dps{\displaystyle}
\def\bve{\;|\;}

\def\<{\langle}
\def\>{\rangle}

\def\H{\mathcal{H}}

\def\xz{x_0}

\def\dim{n}

\newcommand{\pds}[2]{ \langle {#1}\mid#2 \rangle}

\def\JJ{J}
\def\condone{(C)}
\def\condtwo{(C^\star)}

\begin{document}

\title{Gradient flows, second order gradient systems and convexity}

\author{Tahar Z.  Boulmezaoud        \and
        Philippe Cieutat \and Aris Daniilidis}


\institute{Tahar Zamene Boulmezaoud\at
Laboratoire de Math\'ematiques de Versailles, \\
Universit\'e de Versailles Saint-Quentin-en-Yvelines - Universit\'e Paris-Saclay\\
45, avenue des Etats-Unis, 78035, Versailles, Cedex.\\
           \email{\texttt{tahar.boulmezaoud@uvsq.fr}}           
           \and
           Philippe Cieutat\at
Laboratoire de Math\'ematiques de Versailles, \\
Universit\'e de Versailles Saint-Quentin-en-Yvelines - Universit\'e Paris-Saclay\\
45, avenue des Etats-Unis, 78035, Versailles, Cedex. \\
\email{\texttt{philippe.cieutat@uvsq.fr}}
\and 
Aris Daniilidis \at 
DIM--CMM, UMI CNRS 2807, Beauchef 851 (Torre Norte, piso~5)\\
Universidad de Chile, 8370459, Santiago, Chile \\ 
\email{\texttt{arisd@dim.uchile.cl}} \quad Personal page: \texttt{http://www.dim.uchile.cl/{\raise.17ex\hbox{$\scriptstyle\sim$}}arisd} \smallskip \\ 
Research supported by grants BASAL AFB-170001, FONDECYT 1171854 (Chile) and MTM2014-59179-C2-1-P (MINECO of Spain and ERDF of EU).}

\maketitle 

\abstract{We disclose an interesting connection between the gradient flow of a $\C^2$-smooth function $\psi$ and strongly evanescent orbits of the second order gradient system defined by the square-norm of $\nabla\psi$, under adequate convexity assumption. As a consequence, we obtain the following surprising result for two $\C^2$, convex and bounded from  below functions $\psi_1$, $\psi_2$: if $||\nabla\psi_1||=||\nabla\psi_2||$, then $\psi_1=\psi_2 + k$, for some $k\in \R$. }

\tableofcontents


\section{Introduction}
We are interested in the first order gradient system
\begin{equation} \tag{DS-1}\label{Equaf1}
\solf^{\prime}(t)=-\nabla \psi(\solf(t)) , \quad t\geq 0, \; 
\end{equation}
in comparison with the second order gradient system 
\begin{equation} \tag{DS-2}\label{EquaV1}
\solV^{\prime\prime}(t)=  \nabla V(\solV(t)) , \quad t\geq 0, \; 
\end{equation}
where $\psi:\H\to\R$ is a $\C^2$ function (respectively, $V:\H\to\R$ is a $\C^1$ function), $\nabla \psi$, $\nabla V$ denote the respective gradients and $\H$ stands for a Hilbert space, with inner product $\PS{\cdot}{\cdot}$ and associated norm $\norm{\cdot}$. Throughout this work, the functions $\psi$ and $V$ will be linked with the relation
\begin{equation}\label{Vx} 
V(\x) = \frac{1}{2}\norm{\grad \psi(x)}^2, \quad  x \in \H\,. 
\end{equation}
The second-order system \eqref{EquaV1} is introduced here (and studied for the potential $V$ given by \eqref{Vx}) for the first time in the literature. \smallskip\newline
In the sequel the set of critical points of $\psi$ (singular set) will be denoted by
$$
\mathrm{Crit}_{\psi} = \{ \x \in \H \bve \grad \psi(\x) = 0\} = \{ \x \in \H \bve V(\x) = 0\}.  
$$
When $\psi$ is convex, the set $\mathrm{Crit}_{\psi}$ is convex and consists of all (global) minimizers of $\psi$. Therefore, in this case the set of critical values $\psi(\mathrm{Crit}_{\psi})$ is either empty or singleton. We may also observe that $\mathrm{Crit}_{\psi}$ is also the set minimizers of $V$. Therefore it is also convex, whenever $V$ is assumed so. \vskip 1mm
By a {\it global} solution of \eqref{Equaf1} (respectively, \eqref{EquaV1}) we mean a function $\solf\in \C^1([0,+\infty),H)$ (respectively, $\solV\in \C^2([0,+\infty),H)$) satisfying \eqref{Equaf1} (respectively, \eqref{EquaV1}), for all $t\geq0$. In both cases, we impose the initial condition 
\begin{equation} \tag{$\mathrm{I_0}$}\label{initial_cond}
\solf(0) = \solfz \qquad (\,\text{respectively,}\,\,\, \solV(0)=\solfz\,)
\end{equation}
for some given $\solfz \in \H$. This is very common for \eqref{Equaf1} to obtain unique solutions, whereas for \eqref{EquaV1} an additional condition on the initial velocity $v'(0)$ is normally required. We deliberately refrain from doing so, but instead, we require the solutions of \eqref{EquaV1} to be global on $[0, +\infty)$ and to comply with one of the following \textit{asymptotic conditions}, introduced in the following definition.
\begin{definition}[weakly and strongly evanescent solutions] \label{def-1} A global solution $\solV$ of  \eqref{EquaV1} is called
\begin{itemize}
\item  \textit{weakly evanescent} (in short, \textit{w-evanescent}) if it satisfies 
\begin{equation} \tag{w-EV} \label{w-Ev}
\liminf_{t \to +\infty} \norm{\solV'(t)}     = \liminf_{t \to +\infty}  V(\solV(t)) = 0\,, 
 \end{equation}
\item  \textit{strongly evanescent}  if it satisfies
$$\norm{\solV^{\prime}(\cdot)}  \in L^2(0, +\infty) \quad\text{and}\quad V(\solV(\cdot))  \in L^1(0, +\infty) .$$
\end{itemize}
or equivalently
\begin{equation}\label{asymptotic_base}\tag{EV}
\int_{0}^{+\infty} \left( \norm{\solgen'(t)}^2 + V(\solgen(t))\right) \,dt \, < \, +\infty\,.
\end{equation}
\end{definition}

\begin{remark} 
$\mathrm{(i)}$ Conditions \eqref{w-Ev} and \eqref{asymptotic_base} as well as the associated terminology appear to be new in the literature. Both conditions correspond to a kind of boundary condition of the orbit $\solgen(t)$ at infinity. \\
$\mathrm{(ii)}$ Any strongly evanescent solution of  \eqref{EquaV1} is also w-evanescent. 
\end{remark}

It is straightforward to see that any global solution of \eqref{Equaf1} is also solution of \eqref{EquaV1}. However, this solution might fail to satisfy \eqref{asymptotic_base}. To see this, 
let $\dim=1$ and $\psi(x) = -x^2$, for $x \in \R$, and notice that $\solgen(t) = e^{2 t} \xz$ is solution of \eqref{Equaf1} (and consequently of \eqref{EquaV1}), but \eqref{asymptotic_base} fails, since $\solgen \not \in L^2(0, +\infty)$. Conversely, a solution of \eqref{EquaV1} satisfying \eqref{asymptotic_base} and \eqref{initial_cond} might not be solution of \eqref{Equaf1} since the system \eqref{EquaV1}--\eqref{asymptotic_base} does not distinguish between $\psi$ and $-\psi$. 

\smallskip

Let us further consider the following two conditions:
$$
\condone \; \;  \inf_{z \in \H } \norm{\grad \psi(z)}� = 0\qquad \mbox{and}\qquad \condtwo \;\;  \psi \mbox{ is bounded below.}
$$
By Ekeland's Variational Principle \cite[Corollary~2.3]{ekelandVP} we deduce $\condtwo \Longrightarrow \condone.$  This latter condition  $\condone$ is necessary for the existence of w-evanescent solutions of \eqref{EquaV1}. \vskip 1mm

A constant function $\solV = \xs$ is  a w-evanescent solution of \eqref{EquaV1} if and only if $\xs \in \mathrm{Crit}_{\psi}$, while $\mathrm{Crit}_{\psi} \ne \emptyset$ clearly implies $\condone$. If in addition $\psi$ is convex, then  $\condtwo$ is also fulfilled. The example of the following convex $\C^2$ function
\begin{equation}\label{exemple_one}
\psi(x) =
\left\{
\begin{array}{l}
- \ln(1-x)\,, \,\,\mbox { if } x \le 0, \smallskip \\
\phantom{-} \frac{1}{2} x^2 + x\,,\quad \mbox{ if } x \ge 0.
\end{array}
\right. 
\end{equation}
shows that $\condone$ and $\condtwo$ are not equivalent, besides the fact that $\psi$ convex (in this case, only $\condone$ holds). 

\vskip 2mm 

{\it Description of the results.}  First order and second order gradient systems have often been explored independently in
the literature (see, e. g., \cite{brezis1973}, \cite{EkBr2}, \cite{HaleRaugel92}, \cite{HuangSZ1}, \cite{ghoussoub2004}, \cite{chj}, \cite{DLS}, \cite{hj}, \cite{attouch} and references therein). A first innovative aspect of this work is to introduce the particular second order ordinary differential equation \eqref{EquaV1}, for a potential $V(\cdot)$ given by \eqref{Vx}, and shed light on its connection with the first order gradient system \eqref{Equaf1} when either $f$ or $V$ is convex. Exploring this link reveals some unexpected properties of convex functions described below. Another by-product, as we shall see, concerns uniqueness of smooth solutions for certain Eikonal equation.\smallskip \\ 
More precisely,  in this work we show that if either $\psi$ or $V$ is {\it convex}, then any solution of \eqref{EquaV1} satisfying  \eqref{initial_cond}--\eqref{asymptotic_base} is also solution of \eqref{Equaf1}--\eqref{initial_cond}, 
and vice-versa. In particular, the second order system \eqref{EquaV1} coupled with \eqref{initial_cond}-\eqref{asymptotic_base},
is well-posed and can be integrated to obtain the first order system \eqref{Equaf1}.  An important consequence of this result is an intimate link between convexity properties of $\psi$ and of $\norm{\grad \psi}^2$ (Corollary~\ref{VconvFconv}):  
\vskip 1mm
\centerline{($\norm{\grad \psi}^2$  convex and $\psi$ bounded below) $\Longrightarrow$  $\psi$ convex. }  
\vskip 1mm
This leads to the following surprising corollary: 
\vskip 1mm
\centerline{$\norm{\grad \psi_1} = \norm{\grad \psi_2}$ $\Longrightarrow$  $\psi_1 = \psi_2$ +  constant,}
\vskip 1mm
provided that one of the following assumptions is fulfilled:

\begin{itemize}
\item[(a)]  $\psi_1$ and $\psi_2$ are convex  and $\inf \norm{\grad \psi_1} = 0$ (Theorem~\ref{equal_grad_norms_lem}), 
\item[(b)] $\norm{\grad \psi_1}^2$ is convex and $\psi_1$ and $\psi_2$ are bounded below (Corollary~\ref{th_gradgrad_Vcvx}). 
\end{itemize}
\vskip 1mm
Another consequence is a uniqueness property for smooth solutions of the usual 
Eikonal equation 
\begin{equation}\label{eiko_intro}
\|\grad \psi\|^2 = f, 
\end{equation}
in {\it the whole space}. It is well known that uniqueness plays a prominent role in understanding the structure of the set of solutions of \eqref{eiko_intro} (see, e. g., \cite{Kruzkov}, \cite{LionsHJ82}, \cite{CrandallLions1}, \cite{CrandallLions2},  \cite{LPV}, \cite{Barles1994}, \cite{NamahRoq}, 
\cite{BarlesRoq}, \cite{FathiMaderna06},  \cite{ishii09}, \cite{BarlesAL2017} and references therein). Here we obtain uniqueness of bounded below $\C^2$ solutions when $f$ is nonnegative and convex. When  $f$ is only nonnegative, we prove that \eqref{eiko_intro} has most one bounded from below $\C^2$ convex solution. 
If $f$ is only nonnegative, we prove that  any convex and bounded below solution is unique.  
\vskip 1mm

Finally, disclosing the link between \eqref{Equaf1} and \eqref{EquaV1} leads to a simple variational principle for the first order gradient system \eqref{Equaf1}  when
$\norm{\grad \psi}^2$ is convex and $\psi$ bounded below (Proposition~\ref{varia_form_grd}). 

\vskip 2mm

{\it Structure of the manuscript.}  The remainder is organized as follows. In Section~\ref{sec-2} we resume basic properties of the first order system \eqref{Equaf1} for $\psi\in \C^2(\H)$ and for the second order system \eqref{EquaV1} for $V(x)=\frac{1}{2}||\nabla\psi(x)||^2$ that will be used in the sequel. No originality is claimed in Subsection~\ref{ssec-2.1}, as well as in the beginning of Subsection~\ref{ssec-3.1}, where most of the stated properties of the first order system \eqref{Equaf1} are essentially known. These properties are recalled for completeness, provided eventually short proofs to keep the manuscript self-contained. Subsection~\ref{ssec-2.2} contains properties of the system \eqref{EquaV1} with emphasis in Lyapunov functions and in asymptotic behavior of the orbits, while Subsection~\ref{ssec-2.3} is dedicated in comparing the solutions of these two systems.\smallskip

The main results are resumed in Section~\ref{sect_main_res} and organized as follows: Subsection~\ref{ssec-3.1} ensembles all results obtained under the driving assumption that $\psi$ convex, while Subsection~\ref{ssec-3.2} does the same under the assumption $V$ convex. We quote in particular Theorem~\ref{equal_grad_norms_lem} (determination of a convex function by the modulus of its gradient) and its variant Corollary~\ref{th_gradgrad_Vcvx} which are important consequences of Theorem~\ref{th_fcvx_nexist_of_sol} (equivalence of solutions of \eqref{Equaf1} and \eqref{EquaV1} if $\psi$ is convex) and Proposition~\ref{th_Vcvx_nexist_of_sol} respectively. Finally, in Subsection~\ref{ssec-3.3} we associate to the first order system \eqref{Equaf1} an alternative variational principle, which is in the spirit of the results of this work. \smallskip

We assume familiarity with basic properties and characterizations of convex functions. These prerequisites can be found in the classical books \cite{phelps} or \cite{rock}.


\section{Basic properties of first and second  order gradient systems} \label{sec-2}

\subsection{First order gradient system: basic properties}\label{ssec-2.1}

In this subsection we recall for completeness basic properties of solutions of the first order
gradient system \eqref{Equaf1}, which will be used in the sequel. In this subsection 
the functions $\psi\in \C^2(\H)$ and $V(\cdot)$ given in \eqref{Vx}, are not yet assumed to be convex. 

\begin{lemma}[Lyapunov for \eqref{Equaf1}] \label{lemPropGrad1} 
Let $\solf(\cdot)$ be a maximal solution of  \eqref{Equaf1} defined on $[0,T_{\max})$ where $T_{\max}\in(0,+\infty]$. Then, 
\begin{enumerate}[{\rm (i)}]
\item   $\rho(t):= \psi(\solf(t))$ is nonincreasing  on $[0, T)$ and for every $T < T_{\max}$
\begin{equation}\label{L1(i)}
\int_{0}^{T} \norm{\solf^{\prime}(t)}^2 \, dt = \rho(0)-\rho(T)\,;
\end{equation}

\item $\norm{\solf^{\prime}(\cdot)} \in L^2(0, T_{\max})$
  if and only if 
\begin{equation}\label{condition_infsol}
\inf_{0 \le t < T_{\max}} \psi(\solf(t)) > -\infty. 
\end{equation}     
\end{enumerate}
\end{lemma}

\noindent\textit{Proof.} Since $\rho'(t)=\PS{\nabla \psi(\solf(t))}{\solf^{\prime}(t)} = -\norm{\solf^{\prime}(t)}^2=-\norm{\nabla \psi(\solf(t))}^2 \le 0$
we deduce {\rm (i)}. The second assertion follows by taking the limit as $T\rightarrow T_{\max}$.  \hfill $\square$

\bigskip

\begin{remark}[Strict Lyapunov] \label{lem-1} Assuming $\psi\in\C^2(\H)$ yields that both \eqref{Equaf1} and the equation $w'(t)=\nabla\psi(w(t))$ admit unique solutions under a given initial condition. A standard argument now shows that if the initial condition is not a singular point (that is, $\grad \psi(\solf(0)) \ne 0$), then $\grad \psi(\solf(t)) \ne 0$, for every $t>0$ and $\rho$ is strictly decreasing.  
\end{remark} 

\medskip

\begin{lemma} [maximal nonglobal solutions]\label{globalitylem2}
If $\solf(\cdot)$ is a maximal solution of \eqref{Equaf1} which is not global ({\it i.e.} $T_{\max}<+\infty$), then
\begin{equation} \label{condGloba5}
\inf_{0\leq t<T_{\max}} \psi(\solf(t)) = \lim_{t\to T_{\max}} \psi(\solf(t))= -\infty ,
\end{equation}
and 
\begin{equation}\label{condGloba6}
\int_0^{T_{\max}} \norm{\solf^{\prime}(t)}^2 \,dt = +\infty  .
\end{equation}
\end{lemma}

\noindent\textit{Proof.} In view of Lemma~\ref{lemPropGrad1}\,(i) assertions \eqref{condGloba5} and \eqref{condGloba6} are equivalent.  Assume now that \eqref{condGloba6} does not hold. Then the integral 
$$
\int_0^{T_{\max}} \solf^{\prime}(t) dt, 
$$
converges in $\H$ to the element $u(T_{\max})-u_0$, where $u(T_{\max})=\lim_{t\rightarrow T_{\max}}u(t)$. Moreover $\nabla\psi(u(T_{\max}))\neq 0$ ({\it c.f.} Remark~\ref{lem-1}). Considering the Cauchy problem $w^{\prime}(t) = - \nabla \psi(w(t))$ with initial condition $w(T_{\max}) = \solf(T_{\max})$, we deduce that the (presu\-mably maximal) solution 
$\solf(\cdot)$ can be extended to the right on an interval of the form $[0, T_{\max}+\eps)$ for some $\eps >0$, which is a contradiction. \hfill$\square$

\medskip

\begin{corollary}\label{boundedfrombelow}
If $\psi$ is bounded below, then any maximal solution  $\solf(\cdot)$ of \eqref{Equaf1} is global  and 
 $\norm{\solf^{\prime}(\cdot)} \in L^2(0, +\infty)$.
\end{corollary}

\noindent\textit{Proof.}
If $\psi$ is bounded below, then \eqref{condGloba5} cannot be satisfied, and  the  solution $\solf$
is global. Obviously,  \eqref{condition_infsol} is fulfilled yielding $\norm{\solf^{\prime}(\cdot)} \in L^2(0, +\infty)$.
\hfill$\square$

\medskip

\begin{remark}[grad-coercive functions]  A function $\psi\in\C^1(\H)$ is called {\it grad-coercive}  if 
$\| \grad \psi\|$ is bounded on the sublevel sets $[\psi \le \alpha]:=\{x\in\H:\psi(x)\le\alpha\}$, $\alpha \in \psi(\H)$. \smallskip\\
If $\psi$ is grad-coercive then any maximal solution of \eqref{Equaf1} is global. Indeed, let $\solf(\cdot)$ be a maximal solution defined on $[0, T_{\max})$.  Since  $\solf(t) \in [\psi \le \psi(\solf(0))]$, for all $t\in [0,T_{\max})$, the function $\norm{\grad \psi(\solf(\cdot))}$ is bounded on $[0,T_{\max}]$.  Setting $\ds M = \sup_{0\leq t < T_{\max}}  \norm{\grad \psi(\solf(t))}$, we obtain 
$$
\int_0^{T_{\max}} \norm{\solf^{\prime}(t)}^2 dt = \int_0^{T_{\max}} \norm{\grad \psi(\solf(t))}^2 dt \le M  T_{\max} < +\infty,
$$
which contradicts \eqref{condGloba6}. \hfill$\square$ \smallskip 

\noindent Let us observe that  $\psi$ can be grad-coercive without being bounded from below. A simple example
is the identity function $x\mapsto x$ on $\R$. Similarly, a function which is bounded below is not necessarily
grad-coercive, for example the function $x \mapsto \cos(x^2)$. 
 \end{remark}

\begin{remark} [Relation to other domains] Asymptotic behavior of \eqref{Equaf1} has been studied by several authors in the framework of analytic geometry (see \cite{kurdyka}, \cite{sanz}  {\it e.g.}), in relation to convexity (\cite{baillon78}, \cite{DLS}, \cite{DDDL2015}, \cite{MP1991}), to optimization algorithms (\cite{abb}, \cite{apr}, \cite{bbt} {\it e.g.}) and to PDEs (\cite{chj}, \cite{hj} {\it e.g.}). Roughly speaking, good asymptotic behavior requires a strong structural assumption (analyticity or convexity), see \cite{ama} or \cite[p. 12]{palis} for classical counterexamples. 
\end{remark}


\subsection{Second order system: properties of strongly evanescent solutions} \label{ssec-2.2}

In this subsection we emphasize properties of weakly and strongly evanescent solutions of the second order system
 \eqref{EquaV1}, where $\psi\in\C^2(\H)$  and $$V(x)=\frac{1}{2}||\nabla\psi(x)||^2.$$
 
\begin{lemma}[equality of modula]\label{cie2}
Let $\solV(\cdot)$ be a w-evanescent solution of  \eqref{EquaV1}. Then
\begin{equation}\label{nec_exist_sec_grad}
\norm{\solV'(t)} = \norm{\grad \psi (\solV(t))}, \quad \text{for all }\,t\geq 0\,. 
\end{equation} 
\end{lemma}

\noindent\textit{Proof.} It is easily seen that  $I(t):=\frac{1}{2}\norm{\solV^\prime(t)}^2 - V(\solV(t))$ is a first integral of the system \eqref{EquaV1}, that is, for some $k \in \R$ and all $t\geq 0$ it holds $ \norm{\solV'(t)}^2  = k  +  2V(\solV(t))$.  Taking limit inferior as $t\rightarrow +\infty$ we infer from \eqref{w-Ev} that $k=0$ and the result follows.\hfill$\square$

\bigskip

\begin{lemma}[range of orbits]  If $\mathrm{Crit}_{\psi} = \emptyset$, then the range $\set{\solV(t)}{t\geq 0}$  of any w-evanescent solution $\solV(\cdot)$ 
of \eqref{EquaV1} cannot be relatively compact.
\end{lemma}

\noindent\textit{Proof.}  Let  $\solV(\cdot)$ be a w-evanescent solution of \eqref{EquaV1}. If $\set{\solV(t)}{t\geq 0}$
were relatively compact, then there would exist a sequence $(t_n)_{n \geq 0}$ such $\solV(t_n) \to z_0$ for some
$z_0 \in \H$. By  \eqref{w-Ev} we obtain $V(z_0) = 0$. Therefore $\grad \psi (z_0)=0$, that is, $\mathrm{Crit}_{\psi} \ne \emptyset$, a contradiction. \hfill$\square$

\bigskip

The following proposition assembles properties of the strongly evanescent solutions of \eqref{EquaV1}:

\begin{proposition}[Properties of strongly evanescent solutions]\label{genprop1}
Let $\solV(\cdot)$ be a strongly evanescent solution of \eqref{EquaV1}. Then:
\begin{enumerate}[{\rm (i)}]
\item $\dps{\lim_{t \to +\infty} \psi(\solV(t)) \in\R}$ and 
\begin{equation} 
|\psi(\solV(0)) - \psi(\solV(t))| \le \int_{0}^{t} \norm{\solV'(s)}^2 \,ds\,, \quad \text{for all }\,t\geq 0\,. 
\end{equation}

\item  If $\psi$ is coercive $($i.e., $[ \psi \le \alpha ]$ is bounded, for all $ \alpha \in \psi(\H))$, then  $\solV(\cdot)$  is bounded. \medskip

\item If $\normsb{\nabla^2 \psi(\solV(\cdot))}$ is bounded, then $\ds\lim_{t\to+\infty}\norm{\solV^{\prime}(t)}=\lim_{t\to+\infty}V(\solV(t))=0$.
\medskip

\item The function $$\ds\phi(t):=\solV^{\prime}(t)+ \sigma \nabla \psi(\solV(t)), \quad \sigma \in \{-1, 1\}$$ satisfies $$\ds\phi^{\prime}(t) =  \sigma \nabla^2 \psi(\solV(t)) \phi(t) .$$ \smallskip
\end{enumerate}
\end{proposition}

\noindent\textit{Proof.}  Set $r(t):= \psi(v(t))$, $t\geq 0$. Then $|r'(t)|= \pds{ \solV^\prime(t)}{\grad \psi (\solV(t))}$. By Cauchy-Schwarz inequality and Lemma~\ref{cie2} we get 
$
|r'(t)|\leq\dps{  \norm{\solV^\prime(t)}\norm{\grad \psi(\solV(t))}}  \,=\, \dps{   \norm{\solV^\prime(t)}^2 }.
$
Thus, in view of \eqref{asymptotic_base}, $r' \in L^1(0,+\infty)$ and the limit $\dps{\lim_{t \rightarrow +\infty}r(t) = \lim_{t \rightarrow +\infty} \psi(\solV(t)) }$ exists. Moreover, we have 
$$
|r(t)-r(0)| \le  \int_{0}^{t}  \norm{\solV^\prime(s)}^2 \,ds \, \leq\, \int_{0}^{+\infty}  \norm{\solV^\prime(s)}^2 \,ds\,<\,+\infty\,.
$$
 We easily deduce that the range $\{r(t):\,t\geq 0\}$ is bounded, yielding $\solV(t) \in [ \psi \le \eta]$, for some $\eta>0$ and all  $t \ge 0$. Therefore {\rm (ii)} holds. Differentiating the function $V(x)=\frac{1}{2}||\nabla\psi(x)||^2$ and substituting $x=v(t)$ we deduce
 \begin{equation}\label{8a} 
\norm{\grad V(\solV(t))} \leq  \norm{\nabla^2 \psi(\solV(t))} \norm{\grad \psi(\solV(t))}. 
\end{equation}
On the other hand,  
\begin{equation}\label{8b} 
\dps{ \left | \frac{d}{dt}\left[ V(\solV(t))\right]\right|} =\dps{\left | \PS{\nabla V(\solV(t))}  {\solV^{\prime}(t)}\right | }\, \leq \, \dps{ \|\grad V(\solV(t))\| \, \|\solV'(t)\|} \
\end{equation}
Combining \eqref{8a} with \eqref{8b} and recalling \eqref{nec_exist_sec_grad} and the definition of $V$ we get
\begin{equation}\label{8c} 
\dps{ \left | \frac{d}{dt}\left[ V(\solV(t))\right]\right|} \, \le \, 2  \,  \normsb{\nabla^2 \psi(\solV(t))} \, V(\solV(t)). 
\end{equation}

Since $v(\cdot)$ is strongly evanescent, $V(\solV(\cdot)) \in L^1(0,+\infty)$, while $\normsb{\nabla^2 \psi(\solV(\cdot))}$ is bounded by assumption. We deduce from \eqref{8c} that  $\dps{ \frac{d}{dt} [V(\solV(\cdot))]} \in  L^1(0,+\infty)$. Therefore the limit $\ds \lim_{t \to +\infty} V(\solV(t))$ exists (and necessarily equals zero,
since  $V(\solV(\cdot))\in L^1(0,+\infty)$). Thus {\rm (iii)} holds. Finally, {\rm (iv)} follows from direct calculation, using \eqref{EquaV1} and \eqref{Vx}.
\hfill$\square$

\bigskip
The following proposition will be used in the sequel.

\begin{proposition} [Further asymptotic properties of strongly evanescent solutions]\label{prop1a} Let $v(\cdot)$ be a strongly evanescent solution of \eqref{EquaV1} where $V$ is given by \eqref{Vx}. Then  
\begin{equation} \label{11a}
\frac{||\solV(t)-\solV(0)||}{t}, \,\,\frac{\norm{\solV(t)}}{\sqrt{t^2+1}} \in L^2(0, +\infty) \,\,; \quad  \lim_{t \to +\infty} \frac{\|\solV(t)\|}{\sqrt{t}} = 0
\end{equation}
and for every $t \ge 0$ it holds

\begin{equation} \label{consq_hardy}
\int_{0}^{t} \frac{\|\solV(t)-\solV(0)\|^2}{t^2} \,dt \,\,\le \,\,4 \int_{0}^{t} \|\solV^\prime(t)\|^2 \,dt\,.
\end{equation}
\end{proposition}

\noindent\textit{Proof.}  {\rm (i)} Set  $w(t) = \solV(t) - \solV(0)$, $t \ge 0$ (therefore $\dps{ \lim_{t\to 0^+} \frac{w(t)}{t}=v'(0)}$). Integrating by parts and using Cauchy-Schwarz inequality we obtain for every $t > 0$ 
$$
\begin{array}{rcl}
\dps{ \int_{0}^t \frac{\norm{w(s)}^2}{s^2} \,ds } &=& \dps{ - \frac{\norm{w(t)}^2}{t} + 2 \int_{0}^t \frac{\pds{w(s)}{w'(s)}}{s} \,ds } 
\,\le \, \dps{     2 \int_{0}^t \frac{\pds{w(s)}{w'(s)}}{s} \,ds} \medskip  \\
&\le & \dps{ 2 \left(  \int_{0}^t \frac{\norm{w(s)}^2}{s^2} \,ds  \right)^{1/2}  \left( \int_{0}^t \norm{w'(s)}^2  \,ds\right)^{1/2},}
\end{array}
$$
yielding
$$
\int_{0}^{t} \frac{\norm{w(s)}^2}{s^2} \,ds \le 4  \int_{0}^{t} \norm{w'(s)}^2 \,ds =  4  \int_{0}^{t} \norm{\solV'(s)}^2 \,ds\,.
$$
Therefore \eqref{consq_hardy} follows. In particular, since $v(\cdot)$ is strongly evanescent solution, we conclude that $\ds(t^{-1} \norm{w(t)} \in L^2(0, +\infty)$ (hence a fortiori, $\ds(t^2+1)^{-1/2} \norm{w(t)} \in L^2(0, +\infty)$). Since $(t^2+1)^{-1/2} \in L^2(0, +\infty)$, we deduce easily that $(t^2+1)^{-1/2} \norm{\solV(t)} \in L^2(0, +\infty)$. \smallskip

{\rm (ii)}  Fix $t_0 > 0$. Then for all $t >  t_0$ we have
$$
\int_{t_0}^{t} \frac{\|\solV(s)\|^2}{s^2} \,ds = - \frac{\|\solV(t)\|^2}{t} + \frac{\|\solV(t_0)\|^2}{t_0} + 2 \int_{t_0}^t  \frac{\pds{\solV(s)}{\solV'(s)}}{s} \,ds. 
$$
Both integrals in the above expression converge as $t \to +\infty$, yielding that $\dps{\lim_{t \to +\infty}  \frac{\|\solV(t)\|^2}{t} }$ also exists. This limit is necessarily zero since $t^{-1} \norm{\solV(t)} \in  L^2(t_0, +\infty)$. \hfill$\square$

\subsection{Comparison of solutions of \eqref{Equaf1} and \eqref{EquaV1}.}\label{ssec-2.3}

We now focus attention upon comparison between solutions the first order system \eqref{Equaf1}  and evanescent solutions of the second order gradient system \eqref{EquaV1}, where $\psi\in\C^2(\H)$ and $V$ is given by \eqref{Vx}.
\vskip 2mm
 
The following result states that each solution $\solf(\cdot)$ of \eqref{Equaf1}
is also a strongly evanescent solution of  \eqref{EquaV1} unless $\lim_{t \to +\infty} \psi(\solf(t)) = -\infty.$  As underlined in the introduction,
the inverse is more complicated: in general, strongly evanescent solutions of \eqref{EquaV1}  are not necessarily solutions of \eqref{Equaf1}. 
Surprisingly,  under a convexity assumption on either $\psi$ or $V$, strongly evanescent solutions of \eqref{EquaV1} are also solutions of \eqref{Equaf1}. 

\begin{lemma}[Characterization of w-evanescent/strongly evanescent solutions]\label{condition_solf_solV}
Let $\solf(\cdot)$ be  a global solution of  \eqref{Equaf1}. Then, 
\begin{enumerate}[{\rm (i)}]
\item $\solf$  is a global solution of  \eqref{EquaV1}.
\item $\solf$ is a w-evanescent solution of \eqref{EquaV1} if and only if
\begin{equation}\label{condition_grad_bis}
\inf_{t \ge 0} \norm{\grad \psi (\solf(t))}= \inf_{z \in \H} \norm{\grad \psi (z)} = 0. 
\end{equation}   
\item $\solf $ is a  strongly evanescent solution of \eqref{EquaV1} if and only if
\begin{equation}\label{condition_grad_grad}
\inf_{t \ge 0} \psi(\solf(t)) > -\infty. 
\end{equation}     
\end{enumerate}                                                          
\end{lemma}     
             
\noindent\textit{Proof.}  Let $\solf(\cdot)$ be a global solution of \eqref{Equaf1}. This is obviously also a global solution of \eqref{EquaV1} and satisfies 
$\norm{\solf^{\prime}(t)}^2=2V(\solf(t))$. Let us first assume that \eqref{condition_grad_bis} holds. If $\grad \psi(\solf(0)) = 0$, then 
$\solf (t) = \solf(0)$ for all $t \ge 0$ and $\solf(\cdot)$ is trivially w-evanescent. If $\grad \psi(\solf(0)) \ne 0$, then $V(\solf(t)) = \frac{1}{2}\norm{\grad \psi(\solf(t))}^2 \ne 0$ for all $t \ge 0$ ({\it c.f.} Remark~\ref{lem-1}), hence for every $s\ge 0$ it holds 
 $$
   \inf_{t \ge 0} V(\solf(t)) = \inf_{t \ge s}  V(\solf(t))  = 0\,,  
 $$
yielding again that $\solf(\cdot)$ is a w-evanescent solution of \eqref{EquaV1}. The converse is obvious, hence {\rm (i)} is established. \smallskip
 
Assertion {\rm (ii)} follows directly from Lemma~\ref{lemPropGrad1}\,(ii).  \hfill$\square$
\bigskip

Combining Lemma \ref{globalitylem2} and Lemma  \ref{condition_solf_solV} we get 
\begin{corollary}\label{bounded_sol}
Any bounded maximal solution of \eqref{Equaf1} is  a  strongly evanescent solution of  \eqref{EquaV1}. 
\end{corollary}
 Combining Corollary~\ref{boundedfrombelow} with Lemma~\ref{condition_solf_solV}  we obtain the following result.
\begin{proposition}\label{lem_exists_fbounded_from_below}
Let $\psi\in\C^2(\H)$ be bounded from below and $V(x):=\frac{1}{2}||\nabla\psi(x)||^2$. Then for every $\xz \in \H$,  \eqref{EquaV1} has at least one strongly evanescent solution satisfying $\solV(0) = \xz$, which coincides with the unique global solution of the first order system \eqref{Equaf1}.  
\end{proposition}

\section{Main results}\label{sect_main_res}

This section contains the main results of the manuscript, which are presented in three subsections. Before we proceed, let us first recall the following continuous form of the classical Opial's lemma \cite{opial} that will be used in the sequel. (See also \cite[Lemma 17.2.5 (p. 704)]{attouch} for a proof.)
\begin{lemma}[Opial type lemma]\label{ContOpial} Let $S$ be a nonempty subset of a Hilbert space $H$ and $w:[0,+\infty)\to \H$ be a map. Assume that for every $z\in S$, the limit $\ds\lim_{t\to+\infty}\norm{w(t)-z}$ exists and is finite and that all weak sequential limits of $w(\cdot)$, as $t\to+\infty$ belong to $S$. Then $w(t)$ converges weakly to a point of $S$ as $t\to+\infty$ .
\end{lemma}

 \subsection{The case $\psi$ convex.} \label{ssec-3.1} 
 
Throughout this subsection we shall assume that the function $\psi\in\C^2(\H)$ is convex and $V$ is given by \eqref{Vx}. We shall be interested in comparing the 
solutions of \eqref{Equaf1} and \eqref{EquaV1}.  The following result is essentially known (see for instance \cite[Thm 3.1-3.2]{brezis1973} for a proof in a more general context of multivalued evolution equations).

\begin{proposition}[Lyapunov functions for \eqref{Equaf1}] \label{fcvx_exist_grd_syst} 
Let $\psi\in\C^1(\H)$ be convex. Then for every initial condition $\xz \in \H$, the unique maximal solution $\solf(\cdot)$ of  \eqref{Equaf1} satisfying \eqref{initial_cond} is global.  Moreover,

\begin{enumerate}[{\rm (i)}]
\item $\rho(t) = \psi(\solf(t))$ is convex, nonincreasing and
 \begin{equation}\label{fcvx_inf_equa_grd}
 \inf_{t \ge 0}  \psi(\solf(t))  = \lim_{t \to +\infty} \psi(\solf(t))= \inf_{z \in \H} \psi(z). 
 \end{equation}
\item For every $y\in \H$ and $t>0$ it holds $$\dps{ \norm{\solf^{\prime}(t)} \leq \norm{\nabla \psi(y)} + \frac{1}{t}\norm{\solf(0)-y}}.$$ 

\item $t \mapsto  \norm{\solf'(t) } = \norm{ \grad \psi(\solf(t)) } $ is  nonincreasing and
 \begin{equation}\label{decre_V}
 \lim_{t \to + \infty} \norm{\solf^{\prime}(t)} =\inf_{z \in \H} \norm{\nabla \psi(z)}.
 \end{equation}
\item  $\| \solf(\cdot) - \xs \|$ is nonincreasing, for every $\xs \in \mathrm{Crit}_{\psi}$.
\end{enumerate}
\end{proposition}

\begin{remark}[Energy function]\label{rm3.3} Under the assumptions of Proposition~\ref{fcvx_exist_grd_syst}, for every $y \in \H$ we set  
$$\dps{E_y(t) := \frac{1}{2} \norm{\solf(t)-y}^2 + \int_0^t (\psi(\solf(s)) - \psi(y)) \,ds}$$
Since $\psi$ is convex, we deduce $\ds E_y^{\prime}(t)  = \PS{\grad \psi(\solf(t))}{y-\solf(t)}  + \psi(\solf(t))-\psi(y) \le 0,$ that is,
$E_y(\cdot)$ is nonincreasing on $[0,+\infty)$. 
\end{remark}

\noindent The following proposition is well-known. It relates the behaviour of the orbits with the critical points of $\psi$. In the sequel we set
$$ \mathrm{dist}(x, \mathrm{Crit}_{\psi}) :=\inf_{y \in \mathrm{Crit}_{\psi}} \|x-y \|\,.$$

\begin{proposition}\label{prop-p2} Let $\psi\in\C^1(\H)$ be convex and $\solf(\cdot)$ a global solution of \eqref{Equaf1}. 
\begin{enumerate}[{\rm (i)}] 
\item If $\mathrm{Crit}_{\psi} \ne \emptyset$, then $\ds\lim_{t\to +\infty} \norm{\solf'(t)} = 0$ and there exists $\xsp \in \mathrm{Crit}_{\psi}$ such that 
$\solf(t) \underset{t \to +\infty}\rightharpoonup \xsp$ (weakly). Moreover, $\rho_{*}(t):=  \psi(\solf(t))-\psi(\xsp) \in L^1(0,+\infty)$  and 
\begin{equation}\label{19}
 \int_{0}^{+\infty}  \left(\psi(\solf(s))-\min \psi\right) \,ds\,\le \, \frac{1}{2}\,\mathrm{dist}(\solf(0), \mathrm{Crit}_{\psi})^2. 
\end{equation}

\item If $\mathrm{Crit}_{\psi} = \emptyset$, then $\dps{ \lim_{t \to +\infty} \|\solf(t)\| = +\infty}$. \smallskip

\item $\solf(\cdot)$ is bounded if and only if $\mathrm{Crit}_{\psi} \ne \emptyset$. 
\end{enumerate}
\end{proposition}

\noindent\textit{Proof.} The first part of assertion {\rm (i)} follows from \cite[Theorem~17.2.7]{attouch}. Fix now any  $\xsp\in\mathrm{Crit}_{\psi}$. Since $E_{\xsp}(t)\leq E_{\xsp}(0)$ and $\psi(\xsp)=\min\psi$, taking the limit as $t\to +\infty$ we deduce
\begin{equation}\label{19h}
\int_{0}^{+\infty}  \left(\psi(\solf(s))-\min \psi\right) \,ds \, \le \, \frac{1}{2} \|\solf(0)-\xsp \|^2. 
\end{equation}
Taking the infimum in \eqref{19h} for $\xsp\in\mathrm{Crit}_{\psi}$, we obtain \eqref{19}. \smallskip

\noindent Assertion {\rm (ii)} follows from \cite[Corollary~17.2.1]{attouch}, while assertion {\rm (iii)} is a straightforward consequence of the last two assertions. \hfill$\square$ 

\bigskip
The following result will play a key role in the sequel.

\begin{proposition}\label{cvx_sol_sol}
Let $\psi\in\C^2(\H)$ be convex and $V(x)= \frac{1}{2}\norm{\nabla\psi(s)}^2$. Then any w-evanescent solution of \eqref{EquaV1} is also a (global) solution of 
the gradient system \eqref{Equaf1}.  
\end{proposition}

\noindent\textit{Proof.}  Let $\solV(\cdot)$ be a w-evanescent solution of  \eqref{EquaV1} and set  
$\phi(t)= \solV^\prime(t) + \grad \psi(\solV(t))$. Then for all $t\geq 0$, it holds 
$\norm{\phi(t)} \le \norm{\solV^\prime(t)} +  \norm{\grad \psi(\solV(t))}. $
By Lemma~\ref{cie2} we deduce
$$
\norm{\phi(t)}  \leq 2  \norm{\grad \psi(\solV(t))} = 2 \norm{v'(t)} \rightarrow 0
$$
whence $\ds\liminf_{t \to +\infty}  \norm{\phi(t)} = 0$, since $\solV(\cdot)$ is a w-evanescent solution.
We also know that $$\phi'(t) = \grad^2 \psi (\solV(t)) \phi(t)$$ (see Proposition~\ref{genprop1}\,(iv)). Thus, 
$$
\frac{d}{dt}\left( \norm{\phi(t)}^2 \right) = 2 \PS{\phi(t)}{\grad^2 \psi (\solV(t)) \phi(t)} \ge 0,
$$
since $\psi$ is convex. Hence $\norm{\phi}^2$ is increasing. Therefore, since $\ds\liminf_{t \to +\infty}  \norm{\phi(t)} = 0$ we deduce $\phi = 0$, which yields that $\solV(\cdot)$ is solution of the first order gradient system \eqref{Equaf1}. \hfill$\square$

\vskip 3mm

We are now ready to state our main results.

\begin{theorem}[second order gradient system; $\psi$ convex]\label{th_fcvx_nexist_of_sol} 
If $\psi\in\C^2(\H)$ is convex, \eqref{EquaV1} has a w-evanescent solution $\solV(\cdot)$ satisfying \eqref{initial_cond} if and only if  $\condone$ holds. Then $\solV(\cdot)$ is unique and is also the unique solution of the first order system \eqref{Equaf1} that satisfies \eqref{initial_cond}. Moreover, 
\begin{enumerate}[{\rm (i)}]
\item $\solV$ is strongly evanescent if and only if $\psi$ is bounded below\; \smallskip
\item $\solV$ is bounded if and only if $\mathrm{Crit}_{\psi} \ne \emptyset$. 
\end{enumerate}
\end{theorem}

\noindent\textit{Proof.} As already mentioned in the introduction, condition {\condone} is necessary for the existence of a  w-evanescent solution of  \eqref{EquaV1}.
Conversely, suppose that  {\condone} is fulfilled. Then there exists a unique global solution $\solf(\cdot)$ of  \eqref{Equaf1} satisfying $\solf(0)=u_0\in \H$ ({\it c.f.} Proposition~\ref{fcvx_exist_grd_syst}). Condition \eqref{condition_grad_bis} is fulfilled, thanks to \eqref{decre_V} and \condone. Thus,  in view of Lemma~\ref{condition_solf_solV},  $\solf(\cdot)$ is a  w-evanescent solution of \eqref{EquaV1} satisfying \eqref{initial_cond}. Uniqueness is straightforward from Proposition~\ref{cvx_sol_sol}. Indeed, any w-evanescent solution of  \eqref{EquaV1} which satisfies \eqref{initial_cond} is necessarily the unique global solution of  \eqref{Equaf1} under the same initial condition \eqref{initial_cond}. Finally, combining \eqref{condition_grad_grad} with \eqref{fcvx_inf_equa_grd} we deduce that this solution is strongly evanescent if and only if $\psi$ is bounded below. From Proposition~\ref{fcvx_exist_grd_syst}, we also deduce that this solution is bounded if and only if $\mathrm{Crit}_{\psi} \ne \emptyset$. \hfill$\square$

\bigskip

To illustrate Theorem~\ref{th_fcvx_nexist_of_sol}  consider the convex $\C^2$ function $\psi$ given in \eqref{exemple_one}. Recall that $\psi$ satisfies $\condone$ but not $\condtwo$. The first-order system $\solf^\prime(t) = - \psi'(\solf(t))$, $\solf(0) = 0$ has the unique solution $\solf(t) = 1 -\sqrt{1+2 t}$, $t\geq 0$, which is also the unique w-evanescent solution of \eqref{EquaV1} (\textit{c.f.} Theorem~\ref{th_fcvx_nexist_of_sol}).
Clearly this solution is not strongly evanescent ($\psi$ is not bounded from below). \vskip1mm

\noindent An immediate consequence of Theorem~\ref{th_fcvx_nexist_of_sol} and Proposition~\ref{fcvx_exist_grd_syst} is the following result.
\begin{corollary} \label{cor_fcvx_nexist_of_sol} Let $\psi\in\C^2(\H)$ be convex, assume $\condone$ holds and let $\solV(\cdot)$ be a w-evanescent solution of \eqref{EquaV1}. Then $\solV(\cdot)$ satisfies the properties stated in Proposition~\ref{fcvx_exist_grd_syst} and Proposition~\ref{prop-p2}.
\end{corollary}

We are ready to state the following surprising consequence.

\begin{theorem}[determination via modulus of gradient]\label{equal_grad_norms_lem}
Let $\psi_1,\psi_2\in C^2(\H)$ be convex and assume
\begin{itemize}
\item $\norm{\grad \psi_1(z)}= \norm{\grad \psi_2(z)}$ for all $z\in \H$ ; \medskip 
\item  $\ds{\inf_{z \in \H} \norm{\grad \psi_1(z)} = 0}$ (this assumption holds whenever $\psi_1$ or $\psi_2$ is bounded below). 
\end{itemize}

\noindent Then,  $\psi_1= \psi_2 + c$ for some constant $c \in \R$.  
\end{theorem}

\noindent\textit{Proof.} Let $\psi_1$ and $\psi_2$ be two convex functions of class $\C^2$ satisfying 
$\norm{\grad \psi_1} =  \norm{\grad \psi_2}$  and $\dps{ \inf_{z \in \H} \norm{\grad \psi_1(z)} = 0}$. Let $\x \in \H$ an arbitrary point and let $\solV(\cdot)$ be the unique weakly evanescent solution of the system
$$
\solV^{\prime \prime}(t)  = \grad V(\solV(t)) \mbox{ for } t \ge 0, \; \solV(0) = \x, \;
$$
with $V = \frac{1}{2}\|\grad \psi_1\|^2 = \frac{1}{2}\|\grad \psi_2\|^2$ ({\it cf.} Theorem~\ref{th_fcvx_nexist_of_sol}).  
Then  $\solV(\cdot)$ is also solution of the first order systems
$$
\solV^\prime(t) = -\grad \psi_1 (\solV(t)), \quad \solV(0) = \x \,,
$$      
and 
$$
\solV^\prime(t) = -\grad \psi_2 (\solV(t)), \quad \solV(0) = \x.
$$ 
Hence $\grad \psi_1 (x) = \grad \psi_2(x)$. Since $x$ is arbitrary, the result follows. \hfill$\square$

\vskip 3mm

\begin{remark} In \cite{bdp} it has been shown that a continuous (respectively, smooth) convex 1-coercive function can be determined (up to a constant) by knowing its subgradients (respectively, gradients) in specific points of its domain (namely, the ones that correspond to strongly exposed points of the epigraph). Theorem~\ref{equal_grad_norms_lem} asserts that a knowledge of the modulus of the gradient (rather than the gradient itself) suffices to determine a $\C^2$ convex function, provided the function is bounded from below.
\end{remark}

\begin{remark} The assumption that $\psi_1$, $\psi_2$ are bounded from below is important. Consider for instance the example of the functions $\psi_1(x)=x$ and $\psi_2(x)=-x$. \smallskip
\end{remark}

A direct consequence of Theorem~\ref{equal_grad_norms_lem} is the following result concerning uniqueness of convex $\C^2$-smooth solutions of the forthcoming Eikonal equation \eqref{23}.

\begin{corollary}[Eikonal equation - I] 
Let $f\in\C^1(\H)$ be nonnegative. Then, the eikonal equation
\begin{equation}\label{23}
\|\grad \psi\|^2 = f,
\end{equation}
has at most one (up to a constant) convex, bounded below solution in $\C^2(\H)$.  
\end{corollary}

\begin{remark}
The above result might appear to be restrictive at a first sight. Indeed, solving \eqref{23} in a viscosity sense leads to the existence of possibly nonsmooth solutions. In particular, if $\H=\R^d$ and $f(x)\geq\alpha >0$, for all $x\in \R^d$, then any viscosity solution of \eqref{23} is unbounded from below (see \cite[Theorem~1.1]{BarlesAL2017} {\it e.g.}). Nonetheless, the case where $f$ is nonnegative and vanishes is actually of big interest for establishing some weak KAM theorems or existence of solutions for ergodic problems associated with first-order Hamilton-Jacobi equations. It is also known that \eqref{23} may have essentially different solutions, see \cite{LionsHJ82} or \cite{LPV}. See also \cite{NamahRoq} and references therein for the periodic case, and \cite{BarlesRoq}, \cite{FathiMaderna06}, \cite{BarlesAL2017} for the unbounded case. 
In the above framework the set of solutions of \eqref{23} is a challenging issue. The above result as well as forthcoming Corollary~\ref{cor-Eik-II} could eventually shed new light in this intriguing issue. 
\end{remark}

Before we finish this section, let us observe that the assumption $\psi_{1},\psi_{2}\in \C^{2}(\H)$ in Theorem~\ref{equal_grad_norms_lem}, although required in this approach (in view of \eqref{EquaV1}), it does not seem to be indispensable for the validity of the result. Indeed the conclusion of Theorem~\ref{equal_grad_norms_lem} seems plausible also for $C^{1}$-convex functions, or even for (nonsmooth) convex continuous functions, under a different approach. We propose below the following conjecture which, if true, would generalize Theorem~\ref{equal_grad_norms_lem}:
\begin{conjecture}
Let $\psi_{1},$ $\psi_{2}:\H\rightarrow\mathbb{R}$ be two (finite)
convex functions bounded from below such that
\begin{equation}
\inf_{p\in\partial\psi_{1}(x)}\,||p||\,=\,\inf_{q\in\partial\psi_{2}
(x)}\,||q||\text{, \quad for all }x\in\H. \label{P}
\end{equation}
Then $\psi_{1}=\psi_{2}+c$ for some constant $c>0.$
\end{conjecture}

\noindent\textit{Proof of the conjecture if }$\mathcal{H}=\mathbb{R}$.
Let us denote by $D$ the set of points where both $\psi_{1},$
$\psi_{2}$ are simultaneously differentiable. Then $\psi_{1}^{\prime},$
$\psi_{2}^{\prime}$ are increasing functions on $D$ and $D$ is dense in $\mathbb{R}$. It follows from \eqref{P} that $\psi_{1}^{\prime}(x)=\sigma(x)\psi_{2}^{\prime}(x),$ for all $x\in D,$ where $\sigma
(x)\in\{-1,1\}.$ Our task is to establish that $\sigma\equiv1,$ that is
$\psi_{1}^{\prime}=\psi_{2}^{\prime}$ on $D.$ Then, since $\psi_{1},$
$\psi_{2}$ are locally Lipschitz (hence absolutely continuous),
the conclusion follows.

\smallskip

Notice that (\ref{P}) yields that $\psi_{1},$ $\psi_{2}$ have a common set of
global minimizers. Denote by $S=\arg\min\psi_{1}=\arg\min\psi_{2}$ this set,
If $S=\emptyset,$ then it is easily seen that the set of all subgradients
$\partial\psi_{i}(\mathbb{R})=\cup_{x\in\mathbb{R}}\partial\psi_{i}(x)$ is
either contained in $(-\infty,0)$ or in $(0,+\infty),$ for $i\in\{1,2\}$. If
$\partial\psi_{1}(\mathbb{R})$ is contained in $(-\infty,0)$ and $\partial
\psi_{2}(\mathbb{R})$ is contained in $(0,+\infty),$ then we would have
$\sigma\equiv-1$ and $\psi_{1}^{\prime}=-\psi_{2}^{\prime}$ on $D.$ Taking
into account that $\psi_{1}^{\prime},\psi_{2}^{\prime}$ are increasing, we
deduce that $\psi_{1}^{\prime},\psi_{2}^{\prime}$ are constant, which is
impossible since $S=\emptyset$ and $\psi_{1},$ $\psi_{2}$ are bounded from
below. Therefore both $\partial\psi_{1}(\mathbb{R})$ and $\partial\psi
_{2}(\mathbb{R})$ are contained in the same interval  $(-\infty,0)$ or in
$(0,+\infty)$ and $\sigma\equiv1.$\smallskip\ 

Consider now the case $S\neq\emptyset.$ If $S=\mathbb{R}$, then $\psi_{1},$
$\psi_{2}$ are constant and the result holds trivially, while for any $x>\sup
S$ (respectively, any $x<\inf S$) we should have $\partial\psi_{i}%
(x)\subset(0,+\infty)$ (respectively, $\partial\psi_{i}(x)\subset(-\infty,0)$)
by monotonicity. Therefore again $\sigma\equiv1$ and the conclusion follows.   \hfill$\square$

\subsection{The case $V$ convex. }  \label{ssec-3.2}

In this subsection the driving assumption is the convexity of the function $V(x)=\frac{1}{2}||\nabla\psi(x)||^2$, where $\psi\in\C^2(\H)$. The focus is again the comparison of the 
solutions of the systems \eqref{Equaf1} and \eqref{EquaV1}. \smallskip

The following result reveals a characteristic property of the solutions of \eqref{EquaV1}, which is reminiscent to an analogous property for the orbits of the first order system with convex potential.

\begin{proposition}[Contraction of solutions of \eqref{EquaV1}] \label{Lm_unique} Let  $\psi\in\C^2(\H)$ and assume that $V(x)=\frac{1}{2}||\nabla\psi(x)||^2$ is convex. 
If $\solV_1$ and $\solV_2$ are two strongly evanescent solutions of equation \eqref{EquaV1}, then the function
\[
q(t):=\frac{1}{2}\norm{\solV_1(t)-\solV_2(t)}^2
\]
is convex and nonincreasing on $[0,+\infty)$. In particular if $\solV_1(0)=\solV_2(0)$, then $\solV_1=\solV_2$. 
\end{proposition}
\noindent\textit{Proof.} It suffices to prove that $q$ is convex and nonincreasing. Differentiating twice and evoking monotonicity of 
$\nabla V$ (see \cite[Ch. 2]{phelps} {\it e.g.}) we get
 \[
q^{\prime\prime}(t)=\PS{\solV_1^{\prime\prime}(t)-\solV_2^{\prime\prime}(t)}{\solV_1(t)-\solV_2(t)}+\norm{\solV_1^{\prime}(t)-\solV_2^{\prime}(t)}^2
\]
\[
=\PS{\nabla V(\solV_1(t))-\nabla V(\solV_2(t))}{\solV_1(t)-\solV_2(t)}+\norm{\solV_1^{\prime}(t)-\solV_2^{\prime}(t)}^2 \geq 0,
\]
which yields convexity of $q$. Let us prove that $q$ is decreasing. By Proposition~\ref{prop1a}, we have 
$$
\int_{0}^{\infty} \frac{q(t)}{t^2+1}dt = \frac{1}{2}\int_{0}^{\infty} \frac{\norm{\solV_2(t)-\solV_1(t)}^2 }{t^2+1}dt < +\infty.
$$
Suppose that there exists $t_0 > 0$ such that $q'(t_0) > 0$. Since $q$ is convex, we would have
$$
q(t) \ge q'(t_0) (t-t_0) + q(t_0), \quad \text{for all } t \ge t_0\,,$$
yielding 
$$
\int_{0}^{\infty} \frac{q(t)}{t^2+1}dt = +\infty\,, \quad \text{a contradiction.}
$$
Hence, $q$ is decreasing and the result follows. \hfill$\square$

\bigskip

\begin{lemma}\label{th_singf_Vcvx} Let $\psi\in\C^2(\H)$, assume $V(x)=\frac{1}{2}\norm{\nabla\psi(x)}^2$ is convex and let $\solV(\cdot)$ be a strongly evanescent solution of \eqref{EquaV1}.
If $\mathrm{Crit}_{\psi} \ne \emptyset$, then 
\begin{enumerate}[{\rm (i)}]
\item  $h(t):= \|\solV(t) - \xs\|$ is nonincreasing, for every $\xs \in \mathrm{Crit}_{\psi}\,;$ \smallskip
\item $\solV(\cdot)$  is bounded\,; \smallskip
\item There exists $\xsp \in \mathrm{Crit}_{\psi}$  such that 
$\solV(t) \underset{t\to +\infty}\rightharpoonup \xsp$ (weakly).  
 \end{enumerate}
\end{lemma}

\noindent\textit{Proof.} Let $\solV(\cdot)$ be a strongly evanescent solution of the system \eqref{EquaV1} and pick any $\xs \in \mathrm{Crit}_{\psi}$. 
Applying Proposition~\ref{Lm_unique} for $u_1(t) = \solV(t)$ and $u_2(t) = \xs$ for $t \ge 0$, we get (i). Since $\mathrm{Crit}_{\psi} \ne \emptyset$,
 (ii) is a direct consequence of (i). Finally, (iii) can be proved in a similar way as in Proposition~\ref{prop-p2}, using Lemma~\ref{ContOpial} and convexity of $V$. The details are left to the reader. \hfill$\square$

\bigskip

\begin{proposition}[second order gradient system; $V$ convex]\label{th_Vcvx_nexist_of_sol}
Let us assume that $V(x)=\frac{1}{2}\norm{\grad \psi(\x)}^2$ be convex and $\psi\in\C^2(\H)$ be bounded from below. Then \eqref{EquaV1} has a unique strongly evanescent solution satisfying \eqref{initial_cond}, which is also the unique solution of \eqref{Equaf1} that satisfies \eqref{initial_cond}. 
\end{proposition}

\noindent\textit{Proof.} From Corollary~\ref{boundedfrombelow}  and Cauchy-Lipschitz
there exists a unique global solution of \eqref{Equaf1} satisfying the initial condition 
\eqref{initial_cond}. According to Lemma~\ref{condition_solf_solV} this solution is also a strongly evanescent solution of 
\eqref{EquaV1}. Uniqueness follows from Proposition~\ref{Lm_unique}.  \hfill$\square$

\bigskip
We obtain the following consequence.

\begin{corollary}[convexity criterium]\label{VconvFconv}
Let $V(x)= \frac{1}{2}\|\grad \psi(\x)\|^2$ be convex and $\psi\in\C^2(\H)$ be bounded below. Then, $\psi$ is convex. 
\end{corollary} 

\noindent\textit{Proof.} Fix $z_1, z_2 \in\H$ and denote by $\solf_1(\cdot)$ and $\solf_2(\cdot)$ solutions of \eqref{Equaf1} with $z_1$ and $z_2$ as initial data. 
Since $\solf_1$ and $\solf_2$ are also strongly evanescent solutions of  \eqref{EquaV1}, we know that the function 
\[
q(t)=\frac{1}{2}\norm{\solf_1(t)-\solf_2(t)}^2 , \quad\text{ for }t\geq0 
\]
is decreasing ({\it c.f.} Proposition~\ref{Lm_unique}). Thus, 
\[
0 \ge q^{\prime}(t)= - \PS{\nabla \psi(\solf_1(t))-\nabla \psi(\solf_2(t))}{\solf_1(t)-\solf_2(t)},
\]
or equivalently, 
\begin{equation*}
\PS{\nabla \psi(\solf_1(t))-\nabla \psi(\solf_2(t))}{{\solf_1(t)-\solf_2(t)}}\geq0\,.
\end{equation*}
Taking the limit as  $t\to 0$ we deduce that 
\begin{equation*}
\PS{\nabla \psi(z_1)-\nabla \psi(z_2)}{z_1-z_2} \geq 0,
\end{equation*}
which yields that $\psi$ is convex (see \cite[Ch. 2]{phelps} {\it e.g.}).  \hfill$\square$
\vskip 3 mm

\begin{remark}
Corollary~\ref{VconvFconv} is false if $\psi$ is not supposed bounded below: Indeed, let $\psi(x) = x^3$. Then 
$V(x) = \frac{1}{2} | \psi^{\prime} (x)|^2 =  \frac{9}{2}  x^4$ is convex, but $\psi$ is not. Another two dimensional example is $\psi(x_1,x_2)=x_1^4-x_2^2$.
\end{remark}

\begin{remark}\label{remVconvexBounded}
If $V(x) = \frac{1}{2}\norm{\grad \psi(\x)}^2$ is convex and $\psi\in\C^2(\H)$ is bounded from below, then combining Corollary~\ref{VconvFconv} 
with Theorem~\ref{th_fcvx_nexist_of_sol} we deduce that every strongly evanescent solution $\solf(\cdot)$ of \eqref{EquaV1} satisfies the assertions of Corollary~\ref{cor_fcvx_nexist_of_sol} (since $\psi$ is convex). In particular, $\solf(\cdot)$ is bounded if and only if $\mathrm{Crit}_{\psi} \ne \emptyset$. 
\end{remark}
  
  \vskip 3 mm

The following result is a direct consequence of Theorem~\ref{equal_grad_norms_lem} and Corollary~\ref{VconvFconv}.

\begin{corollary}\label{th_gradgrad_Vcvx}
Let $\psi_1, \psi_2\in \C^2(\H)$ be bounded below satisfying $$\norm{\grad \psi_1(\x)} =\norm{\grad \psi_2(\x)}, \quad \text{for all} \; x \in \H.$$ Then, if
$V(x)=\frac{1}{2} \|\grad \psi_1 (\x)\|^2 \; \left(= \frac{1}{2}\|\grad \psi_2 (\x)\|^2\,\right) $ is convex, we deduce that both $\psi_1$ and $\psi_2$ are convex and equal (up to a constant). 
\end{corollary}

An illustration of Corollary~\ref{th_gradgrad_Vcvx} is given, in case $\psi_1$ and $\psi_2$ are of the quadratic form
$$
\psi_1(\x) = \frac{1}{2} \pds{\x}{A_1 \x},  \quad\text{and}\quad \psi_2(\x)= \frac{1}{2} \pds{\x}{A_2 \x}, 
$$
where $A_i$ is a symmetric linear bounded operator, for $i\in\{1,2\}$. One can quickly check
that  $\psi_i$ is bounded below if and only if $A_i$ is positive semidefinite. In the latter case
identity $\; \norm{\grad \psi_1} =\norm{\grad \psi_2}$ means that 
$\norm{A_1 \x }=  \norm{A_2 \x }$ for all $\x \in \H$, yielding
$A_1^2 = A_2^2$. Thus, $A_1=A_2$ (since $A_1$ and $A_2$ are
positive semidefinite) and $\psi_1 = \psi_2$. This is in accordance 
with Corollary~\ref{th_gradgrad_Vcvx}. This example also shows the importance 
of the assumption that $\psi_1$ and $\psi_2$ are bounded below. Indeed,
if $A_2=-A_1 \ne 0$, then $\; \norm{\grad \psi_1} =\norm{\grad \psi_2}$
and $\psi_1 - \psi_2$ is not constant. 

\vskip 3 mm
A direct consequence of Corollary~\ref{th_gradgrad_Vcvx} is the following result.

\begin{corollary}[Eikonal equation - II]\label{cor-Eik-II}
Let $f\in\C^1(\H)$ be nonnegative and convex. Then, the eikonal equation
\begin{equation}
\|\grad \psi\|^2 = f,
\end{equation}
has at most one bounded below solution in $\C^2(\H)$ up to an additive constant.  In addition, this solution
is convex. 
\end{corollary}

\subsection{An alternative variational principle for \eqref{Equaf1}}     \label{ssec-3.3}

In \cite{EkBr2}--\cite{EkBr1}, Br\'ezis and Ekeland proved the following variational characterization when $\psi$ is a proper, convex
and lower semicontinuous functional defined on a Hilbert space $\H$. In this case \eqref{Equaf1} becomes 
$$u'(t)\underset{a.e}\in-\partial \psi(u(t))\,,\quad t\geq 0.$$
If $\solf(\cdot)$ is an absolutely continuous solution of the above differential inclusion on $[0, T]$  for some $T > 0$, 
with initial condition \eqref{initial_cond}, then $\solf(\cdot)$ is the unique minimizer of the functional  
$$
{\mathscr J} (u)  = \int_{0}^{T} \left(\psi(u(t)) + \psi^\star(-u'(t))\right) \,dt + \frac{1}{2}\norm{u(T)}^2, 
$$
where $\psi^\star$ designates the Legendre conjugate of $\psi$. We also refer to \cite{auchmuty2007} and \cite{ghoussoub2004} for extensions of this variational principle. 
\smallskip

We now present an alternative variational principle for the first order gradient  system \eqref{Equaf1}. The formulation is based on the connection with the second order system \eqref{EquaV1}. This latter can be seen as the Euler-Lagrange equation associated to a conventional functional.  
 More precisely,  for any real number  $T > 0$,   we consider the functional
 $$\JJ(T; w) =   \int_0^T  \left( \frac{1}{2} \norm{w^{\prime}(t)}^2 +  \frac{1}{2} \|\grad \psi(w(t))\|^2\right) \,dt +  \psi(w(T)). $$
We state the following 
\begin{proposition}[Variational formulation]\label{varia_form_grd}
Let $V(x)= \frac{1}{2}\norm{\grad \psi(\x)}^2$ be convex, $\psi\in\C^2(\H)$ bounded below and $T > 0$. Then $\solf  \in \C^0([0, T]; \H)\cap\C^1((0,T); \H)$ is a solution of \eqref{Equaf1} 
on $[0, T]$  if and only if 
\begin{equation}\label{25}
\JJ(T;\solf) \le \JJ(T; w),
\end{equation}
for all $w \in \C^0([0, T]; \H)\cap\C^1((0,T); \H)$ satisfying $w(0)=\solf(0)$. 
\end{proposition}

\noindent\textit{Proof.} In view of Corollary~\ref{VconvFconv}, $\psi$ is convex. Let $\solf(\cdot)$ solution of \eqref{Equaf1} 
on $[0, T]$ and $w \in C^1([0, T],\H)$ such that $w(0) = \solf(0)$. Set $h = w -\solf$. Then, 
\begin{equation*}
\begin{split}
\JJ(T;w)  - \JJ(T;\solf) = \int_{0}^T \left(\PS{\solf^{\prime}(t)}{h'(t)} + \frac{1}{2}  \norm{h'(t)}^2 +  V(\solf(t)+h(t)) -  V(\solf(t)) \right)  \,dt \\
 + \quad \psi(\solf(T)+h(T)) -  \psi(\solf(T))\,.  
\end{split}
\end{equation*}
Using convexity of $\psi$ and $V$ we deduce
 \begin{equation*}
\begin{split}
\JJ(T;w)  - \JJ(T;\solf) \geq \int_{0}^T \left( \PS{\solf^{\prime}(t)}{h^{\prime}(t)} +  \PS{\grad V(\solf(t))}{h(t)} \right)  \,dt  \\
 + \PS{\grad \psi(\solf(T))}{(h(T)}\,,
\end{split}
\end{equation*}
\noindent Integrating by parts yields
\begin{equation*}
\begin{array}{l}
\JJ(T;w)  - \JJ(T;\solf) \, \geq \\
\\
\phantom{\JJ(T;\solf) \geq\qquad}\dps{ \geq\,\int_{0}^T \PS{-  \solf''(t) +  \grad V(\solf(t))}{h(t)}  \,dt + \pds{\solf'(T) + \grad \psi(\solf(T))}{h(T)}} .
\end{array}
\end{equation*}

Since $\solf(\cdot)$ is solution of \eqref{Equaf1}, it is also solution of \eqref{EquaV1}, therefore 
$$
\int_{0}^T \PS{-  \solf''(t) +  \grad V(\solf(t))}{h(t)}  \,dt + \pds{\solf'(T) + \grad \psi(\solf(T))}{h(T)} = 0\,,
$$
yielding $\ds\JJ(T;w)   \geq \JJ(T;\solf)$.
\vskip 3 mm
Conversely, let $\solf \in C^1([0, T],\H)$. Assume that
$\JJ(T;w) \ge  \JJ(T;\solf)$ for all $w \in C^1([0, T],\H)$ such that
$w(0) = \solf(0)$. By a conventional argument, we know that $\solf$ is of class $\C^2$. Moreover, $\solf$ satisfies the Euler-Lagrange equation
$\solf''(t) = \grad V(\solf(t))$ and the {\it transversality} condition $\solf'(T) + \grad \psi(\solf(T)) = 0$.
Set $\phi (t)= \solf'(t) +  \grad \psi(\solf(t))$ for $t \geq 0$. 
We known that $\phi$ is a solution of the linear differential equation $\ds\phi^{\prime}(t) = \grad^2 \psi (\solf(t)) \phi(t)$ (see Proposition~\ref{genprop1}) with $\phi(T) = 0$, then $\phi$ is the trivial solution $\phi = 0$, that is to mean
$\solf$ is solution of \eqref{Equaf1} on $[0, T]$. This ends the proof. \hfill$\square$

We now consider the functional
 $$
 \JJ^\star_\infty(w)=   \int_0^{+\infty} \left( \frac{1}{2} \norm{w^{\prime}(t)}^2 +  \frac{1}{2} \|\grad \psi(w(t))\|^2\right) \,dt.
 $$
We also state the following 
\begin{corollary}\label{varia_form_grd}
Suppose that $\psi \in\C^2(\H)$  is bounded below, $\mathrm{Crit}_{\psi} \ne \emptyset$ and  $[\x \mapsto \norm{\grad \psi(\x)}^2]$ is convex.
Then, $\solf \in  \C^1([0, +\infty); \H)$ is a global solution of \eqref{Equaf1} if and only if 
\begin{equation}
J_{\infty} (\solf) \le J_{\infty}(w), 
\end{equation}
for any bounded function $w \in \CC^1([0, +\infty);\H)$ with $w(0)=\solf(0)$.
\end{corollary}

\noindent\textit{Proof.}  Let $\solf  \in   \C^1([0, +\infty); \H)$  be a global solution  of  \eqref{Equaf1}. For $T > 0$
 and $z \in   \C^1([0, T]; \H)$ we set
$$\JJ^{\star}(T;z)  = \int_{0}^T \left(\norm{z'(t)}^2 + \norm{\grad \psi(z(t))}^2\right) \,dt.$$
Let $w \in \C^1([0, +\infty); \H)$ be a bounded function satisfying $w(0) = \solf(0)$ and set
set $h = w-\solf$.  Following the proof of inequality \eqref{25}, we obtain 
\begin{equation}
\JJ^{\star}(T;w)  \ge \JJ^{\star}(T;\solf) + \pds{\solf'(T)}{h(T)},
\end{equation}
Let us observe that $\solf$ is bounded and  $\lim_{T \to +\infty} \norm{ \solf'(T)} = 0$ (thanks to Proposition~\ref{prop-p2}).  Thus $h = w-\solf$ is also bounded. Taking the limit when $T \to +\infty$ yields 
\begin{equation}\label{inequa_JJTINFT}
\JJ^\star_\infty(w)  \ge \JJ^\star_\infty(\solf).
\end{equation}

\noindent Conversely, suppose that $\solf \in \C^1([0, +\infty); \H)$  satisfying \eqref{inequa_JJTINFT}
for any bounded  function $w \in \C^1([0, +\infty); \H)$  with  $w(0) = \solf(0)$. Let $ \xs \in \mathrm{Crit}_{\psi} \ne \emptyset$ and consider the function
$$
w_0(t) = e^{-t} (\solf(0) - \xs) + \xs. 
$$
Let us denote by $[\xs, \solf(0)] = \set{\theta (\solf(0) - \xs) + \xs}{\theta\in[0,1]}$ the segment between $\xs$ and $\solf(0)$. Obviously $[\xs, \solf(0)]$ is a compact subset of $\H$ and $w_0(t) \in [\xs, \solf(0)]$ for all $t \ge 0$. We deduce
$$
V(w_0(t))  = V(w_0(t)) - V(\xs) \le \sup_{x\in [\xs, \solf(0)]}\norm{\grad V(x)} \norm{w_0(t)-\xs} .
$$ 
It follows that $\ds V(w_0(t))   \le \sup_{x\in [\xs, \solf(0)]}\norm{\grad V(x)} \norm{w_0(0) - \xs}e^{-t}$. Therefore we obtain $\JJ^\star_\infty(w_0) < +\infty$, whence $\JJ^\star_\infty(\solf) < +\infty$. \smallskip\newline
\noindent Consider now an arbitrary real number $T > 0$ and let $h \in \C^1([0, +\infty); \H)$ having a compact support  included in  $[0, T]$. Then,
$$
 \JJ^\star_\infty(\solf + h) -  \JJ^\star_\infty(\solf) =  \JJ^\star(T; \solf + h) -  \JJ^\star(T; \solf). 
$$
Thus, 
$$
 \JJ^\star_T(\solf + h)  \ge  \JJ^\star_T(\solf). 
$$
From the latter we deduce that $\solf$ satisfies Euler-Lagrange equation $\solf''(t) = \grad V(\solf(t))$ 
on $(0, T)$. Since $T > 0$ is arbitrary, 
$\solf$ is a global solution of  \eqref{EquaV1} on $[0, +\infty)$. Since $\JJ^\star_\infty(\solf) < +\infty$,
it  is also a strongly evanescent solution. In view of Proposition~\ref{th_Vcvx_nexist_of_sol}, $\solf$ is also
solution of  \eqref{Equaf1}. \hfill$\square$
\bigskip
\begin{remark}
In the second part of the latter  proof, we can show that $\JJ^\star_\infty(\solf) < +\infty$ in another way. Indeed, one can choose $w_0$ as the unique strongly evanescent solution of  \eqref{EquaV1} which satisfies 
$w_0(0) = \solf(0)$ (existence of $w_0$ is ensured by Proposition \ref{th_Vcvx_nexist_of_sol}). In view of Remark \ref{remVconvexBounded}, we know that $w_0$ is bounded.
\end{remark}

\textbf{Acknowledgment.} This work has been completed during mutual visits of the first and the third author to the University of Chile and the University of Versailles SQY respectively. The authors wish to thank the host institutions for hospitality. \newline The authors thank the two referees for their constructive comments.

\bibliographystyle{plain}

\begin{thebibliography}{99}

\bibitem{ama}\textsc{P.-A. Absil, R. Mahony, B. Andrews,} Convergence of the iterates of descent methods for analytic cost functions, \textit{SIAM J. Optim.} \textbf{ 16 } (2005), 531--547.


\bibitem {abb}\textsc{F. Alvarez, J.  Bolte, O.  Brahic, } Hessian
Riemannian gradient flows in convex programming. \textit{SIAM J. Control
Optim.} \textbf{43} (2004), 477--501.

\bibitem{attouch} \textsc{H. Attouch, G. Buttazzo, G. Michaille}, \textit{Variational Analysis in Sobolev and BV spaces: Applications to PDEs and optimization}, MOS-SIAM Series on Optimization \textbf{17} (Second Edition), Society for Industrial and Applied Mathematics (SIAM), Philadelphia, PA, Mathematical Optimization Society (2014). 

\bibitem{apr}\textsc{H. Attouch, J. Peypouquet, P. Redont,} A dynamical approach to an inertial forward-backward algorithm for convex minimization, \textit{SIAM J. Optim.} \textbf{ 24 } (2014), 232--256.
 
\bibitem{auchmuty2007} \textsc{G. Auchmuty,} Variational principles for finite dimensional initial
value problems. In \textit{Control methods in PDE-dynamical systems}, pp.
45--56, Contemp. Math., 426, Amer. Math. Soc., Providence, RI, 2007.

\bibitem{bdp}\textsc{M. Bachir, A. Daniilidis, J.-P. Penot}, Lower subdifferentiability and Integration, \textit{Set-Valued Analysis} \textbf{10 }(2002), 89--108. 

\bibitem{baillon78}\textsc{J.-B. Baillon} Un exemple concernant le comportement asymptotique de la solution du probl\`{e}me {$du/dt+\partial\varphi(u)\ni0$}, \textit{J. Funct. Anal. }\textbf{28} (1978), 369--376.

\bibitem {Barles1994}\textsc{G. Barles}, Solutions de viscosit\'e des \'equations de {H}amilton-{J}acobi, "Math\'ematiques \& Applications" [Mathematics \& Applications], Vol. 17,  Springer-Verlag, Paris  (1994). 

\bibitem{BarlesRoq}\textsc{G. Barles and J. M. Roquejoffre}, Ergodic type problems and large time behaviour of unbounded solutions of Hamilton-Jacobi equations  \textit{ Communications in Partial Differential Equations} \textbf{31}, No. 7-9, (2006), 1209--1225. 

\bibitem{BarlesAL2017}\textsc{G. Barles, O. Ley, T.-T. Nguyen and T. Phan},
Large time behavior of unbounded solutions of  first-order Hamilton-Jacobi 
equations in $\R^n$. \textit{ArXiv},  Preprint 1709.08387 (2017). 


\bibitem{bbt}\textsc{H. Bauschke, J. Bolte, M. Teboulle,} A descent lemma beyond Lipschitz gradient continuity: first-order methods revisited and applications, \textit{Math. Oper. Res.} \textbf{42} (2017), 330--348.
 
\bibitem{brezis1973}\textsc{H. Brezis,} \textit{Op\'erateurs maximaux monotones et semi-groupes de
contractions dans les espaces de Hilbert}. North-Holland Mathematics Studies,
No. 5. Notas de Matem\'atica (50). American Elsevier Publishing Co., Inc., New
York, 1973. vi+183 pp.

\bibitem{EkBr2}\textsc{H. Brezis, I. Ekeland,}  Un principe variationnel associ\'e \`{a} certaines
\'equations paraboliques. Le cas d\'ependant du temps. (French) \textit{C. R.
Acad. Sci. Paris S\'er.} A-B \textbf{282} (1976), no. 20, Ai, A1197--A1198.

\bibitem{EkBr1}\textsc{H. Brezis, I. Ekeland,} Un principe variationnel associ\'e \`a certaines
\'equations paraboliques. Le cas ind\'ependant du temps. (French) \textit{C. R.
Acad. Sci. Paris S\'er.} A-B \textbf{282} (1976), no. 17, Aii, A971--A974.

\bibitem{chj}\textsc{R. Chill, A. Haraux, M. Jendoubi,} Applications of the \L ojasiewicz-{S}imon  gradient inequality to gradient-like evolution equations,
\textit{Anal. Appl. (Singap.)} \textbf{7} (2009), 351--372.

\bibitem{CrandallLions1}\textsc{M. G. Crandall and P. L. Lions}, Condition d'unicit\'e pour les solutions g\'en\'eralis\'ees des \'equations de {H}amilton-{J}acobi du premier ordre. \textit{Comptes Rendus des S\'eances de l'Acad\'emie des Sciences. S\'erie I. Math\'ematiques} \textbf{292}, No.  3 (1981), 183-186.

\bibitem{CrandallLions2}\textsc{M. G. Crandall and P. L. Lions}, Viscosity solutions of {H}amilton-{J}acobi equations in infinite dimensions. {VII}. {T}he {HJB} equation is not always satisfied. \textit{Journal of Functional Analysis} \textbf{194}, No.  1 (1994), 111--148.

\bibitem{DLS}\textsc{A. Daniilidis, O. Ley, S. Sabourau,}  Asymptotic behaviour of self-contracted planar curves and gradient orbits of convex functions,
\textit{J. Math. Pures Appl.} \textbf{ 94} (2010), 183--199. 

\bibitem {DDDL2015}\textsc{A. Daniilidis, G. David, E. Durand, A. Lemenant,} Rectifiability of self-contracted curves in the Euclidean space
and applications, \textit{J. Geom. Anal.} \textbf{25} (2015), 1211--1239.

\bibitem{ekelandVP}\textsc{I. Ekeland,} On the variational principle. \textit{J. Math. Anal.
Appl.} \textbf{47} (1974), 324--353.

\bibitem{FathiMaderna06}\textsc{A. Fathi and E. Maderna}, Weak {KAM} theorem on non-compact manifolds. \textit{NoDEA Nonlinear Differential Equations Appl.} \textbf{14}, No. 1-2 (2007), 1--27.

\bibitem{ghoussoub2004}\textsc{N. Ghoussoub, L. Tzou,} A variational principle for gradient flows.
\textit{Math. Ann.} \textbf{330} (2004), no. 3, 519--549.

\bibitem{HaleRaugel92}\textsc{J. K. Hale and G. Raugel}, Convergence in gradient-like systems with applications to {PDE}. \textit{Zeitschrift f\"ur Angewandte Mathematik und Physik. ZAMP.} \textbf{43}, No. 1, (1992), 63--124.

\bibitem{hj}\textsc{A. Haraux, M. Jendoubi,} The  \L ojasiewicz gradient inequality in the infinite-dimensional Hilbert space framework, \textit{J. Funct. Anal.}  \textbf{260 } (2011), 2826--2842.


\bibitem {HuangSZ1}\textsc{S.-Z. Huang}, Gradient inequalities. Mathematical Surveys and Monographs, \textbf{126}, American Mathematical Society, Providence, RI (2006).

\bibitem{ishii09}\textsc{H. Ishii} Asymptotic solutions of Hamilton-Jacobi equations for large time and related topics.  \textit{In ICIAM07 6th International Congress on Industrial and Applied Mathematics},  Eur. Math. Soc., Zurich, 2009 (2009), 193-217.


\bibitem{Kruzkov}\textsc{S. N. Kru\v{z}kov}, Generalized solutions of {H}amilton-{J}acobi equations of Eikonal type. I. \textit{Mat. Sb. (N.S.)},\textbf{98}, No. 3   (1975),  450--493. 


\bibitem{kurdyka}\textsc{K. Kurdyka,} On gradients of functions definable in o-minimal structures, \textit{Ann. Inst. Fourier (Grenoble)} \textbf{48 } (1998), 769--783. 

\bibitem{LPV}\textsc{P. L. Lions, G. Papanicolau and  S. R. S. Varadhan}, Homogenization of Hamilton-Jacobi equation. Unpublished preprint (1987). 

\bibitem{LionsHJ82}\textsc{P. L. Lions}, Generalized solutions of {H}amilton-{J}acobi equations, Research Notes in Mathematics, Vol. 69,  Springer-Verlag, Pitman, Boston, Mass.-London (1982). 

\bibitem {MP1991}\textsc{P. Manselli, C. Pucci, } Maximum lenght of steepest curves for quasiconvex functions, \textit{Geometriae Dedicata} \textbf{38}
(1991), 211--227.

\bibitem{NamahRoq}\textsc{G. Namah and J. M. Roquejoffre}, Remarks on the long time behaviour of the solutions of {H}amilton-{J}acobi equations. \textit{Communications in Partial Differential Equations} \textbf{24}, No. 5-6, (1999), 883-893.


\bibitem{opial}\textsc{Z. Opial}, Weal convergence of the sequence of successive approximations for nonexpansive mappings, \textit{Bull. Amer. Math. Soc.} \textbf{73} (1967), 591--597.

\bibitem {palis}\textsc{J. Palis, W. De Melo,} \textit{Geometric theory of
dynamical systems. An introduction,} (Translated from the Portuguese by A. K.
Manning), Springer-Verlag, New York-Berlin, 1982.


\bibitem{phelps} \textsc{R. Phelps,}  \textit{Convex functions, monotone operators and differentiability.} (2nd eds). Lecture Notes in Mathematics \textbf{1364}, Springer-Verlag, Berlin, 1993.

\bibitem{rock}\textsc{T. Rockafellar,} \textit{Convex analysis.} (Reprint of the 1970 original.) Princeton University Press, Princeton, NJ, 1997.

\bibitem{sanz}\textsc{F. Sanz},  Non-oscillating solutions of analytic gradient vector fields, \textit{Ann. Inst. Fourier (Grenoble)} \textbf{ 48} (1998), 1045--1067.


\end{thebibliography}

\end{document}